\documentclass[12pt]{article}
\title{On coefficients of Yablonskii-Vorob'ev polynomials}
\author{Masanobu Kaneko and Hiroyuki Ochiai}
\date{}
\newcommand{\address}{\noindent
Email: mkaneko@math.kyushu.u-ac.jp \\
Graduate School of Mathematics, Kyushu University, \\
Fukuoka 812-8581, Japan. \\
\\
E-mail: ochiai@math.titech.ac.jp \\
Department of Mathematics, Tokyo Institute of Technology, \\
Meguro, Tokyo 152-8551, Japan.}

\usepackage{amsfonts}
\usepackage{amsmath}
\usepackage{amssymb}
\usepackage{theorem}
\newtheorem{theorem}{{\bf Theorem}}
\newtheorem{proposition}[theorem]{{\bf Proposition}}
\newtheorem{corollary}[theorem]{{\bf Corollary}}
\newtheorem{lemma}[theorem]{{\bf Lemma}}

\theorembodyfont{\rmfamily}
\newtheorem{remark}[theorem]{{\bf Remark}}
\newtheorem{example}[theorem]{{\bf Example}}
\newcommand{\qed}{\hspace*{\fill} \fbox{\;}\par\bigskip}

\newcommand{\Z}{{\mathbf Z}}
\newcommand{\Q}{{\mathbf Q}}
\newcommand{\F}{{\mathbf F}}
\newcommand{\p}{{\mathfrak{p}}}
\newcommand{\Zp}{{\Z_{(p)}}}
\newcommand{\barh}{{\widetilde{h}}}

\newcommand{\ta}{{\widetilde a}}
\newcommand{\tb}{{\widetilde b}}
\newcommand{\tc}{{\widetilde c}}

\begin{document}
\maketitle
{\bf Abstract.}
We give a formula for the coefficients of the Yablonskii-Vorob'ev 
polynomial. Also the reduction modulo a prime number of the 
polynomial is studied.%
\footnote{
{\it 2000 Mathematical Subject Classification}. 
Primary 34M55 ; Secondary 33E17.

{\it Key Words and Phrases}.
Yablonskii-Vorob'ev polynomial, Schur function,
Painlev\'e equation.

%

}

%
\section{Introduction}

The object of study in the present article is a sequence of polynomials
$T_n(x)\in\Z[x]$ ($n=0,1,2,\ldots$), referred to as the
Yablonskii-Vorob'ev polynomials,  satisfying the recursion
\begin{equation}\label{eq:recursion}
T_{n+1}(x) T_{n-1}(x) = x T_n(x)^2 +T_n(x) T''_n(x) - {T'_n(x)}^2,
\end{equation}
with the initial condition $T_0(x)=1$, $T_1(x)=x$. The first few are
\begin{eqnarray*}
  T_2(x)&=&x^3-1,\\
T_3(x)&=&x^6-5x^3-5,\\
T_4(x)&=&x^{10}-15x^7-175x,\\
T_5(x)&=&x^{15}-35x^{12}+175x^9-1225x^6-12250x^3+6125.
\end{eqnarray*}
Note that we have adopted a normalization  different from the usual
one (see the remark at the end of Section 2). 

Although it is not clear
a priori that the recursion \eqref{eq:recursion}  gives a sequence of
{\it polynomials}, we know it does indeed, the fact which is  most
naturally  explained in the context of connection with rational
solutions of the  second Painlev\'e equation ($P_{\scriptsize{II}}$).
(See, e.g., \cite{FOU,Taneda} for this and related subjects.)
Specifically, the logarithmic derivative
$y=T'_{n}(x)/T_{n}(x)-T'_{n-1}(x)/T_{n-1}(x)$ of the ratio $T_{n}(x) /
T_{n-1}(x)$ is a solution of
$$
\frac{d^2y}{dx^2}=2y^3-4xy+4n. \leqno{(P_{\scriptsize{II}})} $$
As such, the Yablonskii-Vorob'ev polynomial can be thought of as a
non-linear analogue of the classical special polynomials associated to
linear differential equations. In this paper, we  discuss some
properties including explicit formulas and reductions modulo primes of
coefficients of this ``Painlev\'e special polynomial''. We note that,
owing to the connection with Schur functions, such results also give
a kind of information on certain character values of irreducible
representations of symmetric groups.

Now we state our main results. Using the recursion
\eqref{eq:recursion}, it is easy to see by induction that the
polynomial $T_n(x)$ is monic of degree $n(n+1)/2$ and has the
following  expansion;
\begin{equation}
  \label{eq:tn}
  T_n(x)=\sum_{j\ge0} t_j(n)x^{3j+\delta},\quad t_j(n)\in\Z,
\end{equation}
where $\delta=1$ if $n\equiv 1\bmod 3$ and $0$ otherwise.
Set $$\mu_n = \prod_{k=1}^{n} (2k-1)!!.$$  The first theorem gives
the coefficient of the term of lowest degree ($=$ constant term if
$n\equiv0,2\bmod3$ and the term of degree $1$ if $n\equiv1\bmod3$) of
$T_n(x)$.
\begin{theorem}\label{th:1} We have
\begin{eqnarray}
\label{const}
t_0(n) =\left\{ \begin{array}{ll}
(-1)^m 3^{-(3m-1)m/2} \mu_n/(\mu_{m-1}^2\mu_m), & \quad
\mbox{if $n=3m-1$},\\
(-1)^m 3^{-(3m+1)m/2} \mu_n/(\mu_{m-1}\mu_m^2), & \quad
\mbox{if $n=3m$},\\
(-1)^m 3^{-3(m+1)m/2} \mu_n/\mu_m^3, &\quad
\mbox{if $n=3m+1$}.
\end{array} \right.
\end{eqnarray} 
\end{theorem}
As for the higher coefficients, we show the following.
\begin{theorem}\label{th:2}
For fixed $j$, the function $n\mapsto t_j(n)/t_0(n)$ extends to
a polynomial function in $n$.
\end{theorem}
Several examples of the theorem will be given at the end of Section 3.

The next result concerns the reduction modulo a prime of 
the polynomial $T_n(x)$.
\begin{theorem}\label{th:mp+n}
  For a prime number $p>3$ and any non-negative integers $m$ and $n$,
  we have
\[
T_{mp+n}(x) \equiv x^{d_{mp+n}-d_n} T_n(x) \bmod p,
\]
where $d_n=n(n+1)/2$, the degree of $T_n(x)$.
\end{theorem}

\section{Constant terms}

To prove Theorem~\ref{th:1}, we recall the determinant expression of
the Yablonskii-Vorob'ev polynomial of Jacobi-Trudi type
\cite{Kajiwara}.  Define a family of polynomials $h_k(x) \in \Q[x]$
($n=0,1,2,\ldots$) by the generating function
\begin{equation}\label{eq:generating}
e^{x \lambda + \frac13 \lambda^3} 
= \sum_{k=0}^\infty h_k(x) \lambda^k,
\end{equation}
and set $h_{-1}=h_{-2}= \cdots =0$. Writing the left-hand side 
as $e^{x\lambda}e^{\lambda^3/3}$ and expanding this out, we see that
the polynomial $h_k(x)$ is given by 
\begin{equation}\label{eq:hk}
 h_k(x) = \sum_{i=0}^{[k/3]} \frac{1}{3^ii! (k-3i)!} x^{k-3i},
\end{equation}
where $[k/3]$ is the greatest integer which does not exceed $k/3$.
In particular, the degree of $h_k(x)$ is $k$ and the leading coefficient
is $1/k!$.  Set
\begin{equation}\label{eq:schur}
\tau_n(x) = \det( h_{j-2i+n+1}(x) )_{1 \le i,j \le n}.
\end{equation}
The polynomial $\tau_n(x)$ is known as the 2-core Schur polynomial
attached to the staircase partition of depth $n$.  The degree of
$\tau_n(x)$ is at most $d_n=n(n+1)/2$ since the degree of $h_k(x)$ is $k$,
but it turns out that it is exactly $d_n$ and the coefficient of 
$x^{d_n}$ in $\tau_n(x)$ is given by $\mu_n^{-1}=1/\prod_{k=1}^n(2k-1)!!$,
as the following lemma shows.
\begin{lemma}\label{lemma:highest} We have
$$\det(1/(j-2i+n+1)!)_{1 \le i,j \le n} = \mu_n^{-1},$$
where we understand $1/l!=0$ if $l<0$.
\end{lemma}
A proof is found in \cite{Stanley}, Corollary 7.16.3 (formula 7.71)
combined with Corollary 7.21.6.

The determinant formula for the Yablonskii-Vorob'ev polynomial
\cite{Kajiwara, Yamada} asserts that $T_n(x)$ is a constant multiple of
$\tau_n(x)$:
\begin{equation}\label{eq:detexp}
T_n(x)=\mu_n \tau_n(x).
\end{equation}
\noindent{\it Proof of Theorem~\protect{\ref{th:1}}}.

Suppose $n=3m-1$. Then $t_0(n)$ is the constant term of $T_n(x)$.
 From equations \eqref{eq:detexp} and \eqref{eq:schur}, we want to
compute the determinant
$$\tau_n(0) = \det(h_{j-2i+3m}(0))_{1 \le i,j \le 3m-1}.$$
The point is that this determinant splits into three blocks and we can 
calculate each block separately by using Lemma~\ref{lemma:highest}.  Actually,
noting from \eqref{eq:hk} that $h_{3l}(0)=1/(3^ll!)$ and
$h_{3l-1}(0)=h_{3l+1}(0)=0$, we proceed as follows:

\begin{itemize}
\item[(1)] For $i=3k$ with $1 \le k \le m-1$, the $(i,j)$ entry
  $h_{j-6k+3m}(0)$ is zero unless $j=3l$ with $1\le l \le m-1$,
  in which case the value is $h_{3(l-2k+m)}(0) =
  1/(3^{l-2k+m}(l-2k+m)!)$.  Then, by Lemma~\ref{lemma:highest}, the
  determinant of $m-1$ by $m-1$ matrix with these $(k,l)$ entries is
  equal to $1/(3^{(m-1)m/2}\mu_{m-1})$.
  
\item[(2)] For $i=3k-1$ with $1 \le k \le m$, the $(i,j)$ entry
  $h_{j-6k+3m+2}(0)$ is zero unless $j=3l-2$ with $1\le l \le
  m$, in which case the value is $h_{3(l-2k+m)}(0) =
  1/(3^{l-2k+m}(l-2k+m)!)$.  Noting that this is equal to $0$  for
  $k=m$ and $l<m$, and $1$ for $k=l=m$, we see that the $m$ by $m$
  determinant is equal to the one in (1) as above, i.e., equal to
  $1/(3^{(m-1)m/2}\mu_{m-1})$.
  
\item[(3)] Similarly, for $i=3k-2$ with $1 \le k \le m$, the $(i,j)$
  entry $h_{j-6k+3m+4}(0)$ is zero unless $j=3l-1$ with $1\le l
  \le m$.  By Lemma~\ref{lemma:highest}, the determinant of $m$ by $m$
  matrix with entries $1/(3^{l-2k+m+1}(l-2k+m+1)!)$ is equal to
  $1/(3^{(m+1)m/2}\mu_{m})$.
\end{itemize}
Combining the above three, we conclude $\tau_{3m-1}(0)=  \pm1/
(3^{(3m-1)m/2}\mu_{m-1}^2\mu_{m})$, the sign being  the inversion
number of the permutations of rows and columns, which, as is readily
seen, is equal to $(-1)^m$. This establishes the formula in the case
of $n=3m-1$.  The computation in the case when $n=3m$ is exactly the
same and will be omitted.  When $n=3m+1$, $t_0(n)$ is not the constant
term and the above  computation does not work. But the following lemma
allows us to reduce this case to the preceding two.
\begin{lemma}\label{lemma:wronskian} We have
  $$T_{n-1}(x)T_{n+1}'(x)-T_{n-1}'(x)T_{n+1}(x)=(2n+1)T_n(x)^2$$
  for
  all $n$.
\end{lemma}
See \cite[p.92]{Taneda} or \cite[p.188]{FOU} for a proof. Putting $n=3m$
in the lemma and comparing the constant term of both sides, we obtain
\begin{equation}
  \label{eq:key}
  t_0(3m-1)t_0(3m+1) =(6m+1)t_0(3m)^2.
\end{equation}
 From this, we have
\begin{eqnarray*}
  t_0(3m+1)&=&(6m+1)t_0(3m)^2/t_0(3m-1)\\
&=&(-1)^m(6m+1)3^{-3(m+1)m/2}\mu_{3m}^2/(\mu_{3m-1}\mu_m^3)\\
&=&(-1)^m3^{-3(m+1)m/2}\mu_{3m+1}/\mu_m^3,
\end{eqnarray*}
which completes the proof of Theorem~\ref{th:1}. \qed

\begin{remark}
  When $n\equiv 0,2\bmod 3$, there is an alternative way to derive the
  formula in Theorem~\ref{th:1} from the hook-type formula of $T_n(x)$
  in \cite{Yamada} (the authors would like to thank Masatoshi Noumi
  for pointing out this). However, the case $n\equiv 1\bmod 3$ does
  not follow from the hook-type formula.
\end{remark}
\begin{remark}\label{remark:3}
  As mentioned in the introduction, the usual recursion for the
  Yablonskii-Vorob'ev polynomials is
\begin{equation}\label{eq:recursion'}
T_{n+1}(x) T_{n-1}(x) = x T_n(x)^2 -4(T_n(x) T''_n(x) - {T'_n(x)}^2).
\end{equation}
If in general we start with the recursion
\begin{equation}\label{eq:gen-rec}
T_{n+1}(x) T_{n-1}(x) = x T_n(x)^2 +a(T_n(x) T''_n(X) - {T'_n(x)}^2),
\end{equation}
$a$ being a constant, and the same initial values $T_0(x)=1$ and $T_1(x)=x$,
the formula for the lowest term in Theorem~\ref{th:1} changes only by
a power of $a$, namely,
\begin{eqnarray*}
t_0(n) =\left\{ \begin{array}{ll}
(-1)^m (a/3)^{(3m-1)m/2} \mu_n/( \mu_{m-1}^2\mu_m), & \quad
\mbox{if $n=3m-1$},\\
(-1)^m (a/3)^{(3m+1)m/2} \mu_n/( \mu_{m-1}\mu_m^2),&\quad
\mbox{if $n=3m$},\\
(-1)^m (a/3)^{3(m+1)m/2} \mu_n/\mu_m^3,&\quad
\mbox{if $n=3m+1$}.
\end{array} \right.
\end{eqnarray*}
\end{remark}

\section{Higher coefficients}

For the proof of Theorem~\ref{th:2}, it is convenient to use
different symbols for $t_j(n)$ according to the congruence classes
of $n$ modulo $3$. Put $$a_j(m)=t_j(3m-1),\  b_j(m)=t_j(3m),\ \ 
\mbox{and}\ \ c_j(m)=t_j(3m+1).$$ 
Also put $$\ta_j(m)=a_j(m)/a_0(m),\  \tb_j(m)=b_j(m)/b_0(m),
\ \ \mbox{and}\ \  \tc_j(m)=c_j(m)/c_0(m).$$ 

\noindent{\it Proof of Theorem~\ref{th:2}}.
First let $n=3m$. We substitute the expansion \eqref{eq:tn} 
into the recursion \eqref{eq:recursion} and compare the coefficients
of $x^{3k+1}$ for $k \ge 0$ to obtain
\begin{eqnarray*}
\sum_{i=0}^k c_i(m) a_{k-i}(m)
&=& \sum_{i=0}^k b_i(m) b_{k-i}(m) + \sum_{i=1}^{k+1} 3i(3i-1) b_i(m) 
b_{k+1-i}(m)\\
&-& \sum_{i=1}^k 9 i j b_i(m) b_{k+1-i}(m).
\end{eqnarray*}
Dividing both sides by $c_0(m)a_0(m)$, which is equal to
$(6m+1)b_0(m)^2$ by \eqref{eq:key}, and separating the term with
$i=k+1$ in the middle sum on the right (the only place where $b_{k+1}(m)$
appears), we obtain
\begin{eqnarray}\label{eq:ind-B}
3(k+1)(3k+2) \tb_{k+1}(m) 
&=& (6m+1) \sum_{i=0}^k \tc_i(m) \ta_{k-i}(m)
- \sum_{i=0}^k \tb_i(m) \tb_{k-i}(m)\nonumber \\ 
&&+3\sum_{i=1}^k  i(3k-6i+4) \tb_i(m)\tb_{k+1-i}(m)
\end{eqnarray}
for $k\ge0$.
Similarly, for $n=3m-1$ we obtain from the recursion
\eqref{eq:recursion}
\begin{eqnarray}\label{eq:ind-A}
&& 3(k+1)(3k+2) \ta_{k+1}(m)\\ 
&=& -(6m-1) \sum_{i=0}^k \tc_i(m-1) \tb_{k-i}(m)
- \sum_{i=0}^k \ta_i(m) \ta_{k-i}(m)\nonumber \\ 
&&+3\sum_{i=1}^k  i(3k-6i+4) \ta_i(m)\ta_{k+1-i}(m) \nonumber
\end{eqnarray}
for $k\ge0$.  Here, we have used the identity
$b_0(m)c_0(m-1)=-(6m-1)a_0(m)^2$ which follows from
Lemma~\ref{lemma:wronskian} by putting $n=3m-1$ and comparing the
constant terms of both sides.  For $n=3m+1$, we compare the constant
terms in the recursion \eqref{eq:recursion} to get
$a_0(m+1)b_0(m)=-c_0(m)^2$, and then with this we obtain as above
(comparing the coefficients of $x^{3k+3}$ in \eqref{eq:recursion})
\begin{eqnarray}\label{eq:ind-C}
&&(3k+1)(3k+4) \tc_{k+1}(m) \\
&=& -\sum_{i=0}^{k+1} \ta_i(m+1) \tb_{k+1-i}(m)
- \sum_{i=0}^{k} \tc_i(m) \tc_{k-i}(m)\nonumber \\ 
&&+\sum_{i=1}^{k}  (3i+1)(3k-6i+4) \tc_i(m)\tc_{k+1-i}(m)\nonumber
\end{eqnarray}
for $k\ge 0$. 

Now we prove Theorem~\ref{th:2} by induction on $j$. For $j=0$, the
required property, which is equivalent to the statement that
$\ta_j(m), \tb_j(m)$, and  $\tc_j(m)$ are polynomials in $m$,
holds trivially.  Suppose the property holds up to $j\le k$.  Then
equations \eqref{eq:ind-A} and \eqref{eq:ind-B} ensures respectively
that both $\ta_{k+1}(m)$ and $\tb_{k+1}(m)$ are polynomials in $m$.
Then, we conclude in turn by equation (\ref{eq:ind-C}) that 
$\tc_{k+1}(m)$ is also a polynomial in $m$. This completes the proof 
of Theorem~\ref{th:2}.  \qed

Equations \eqref{eq:ind-B}, \eqref{eq:ind-A}, and \eqref{eq:ind-C}
allows us to compute explicitly the polynomials $\ta_j(m), \tb_j(m)$,
and $\tc_j(m)$. First several examples are given below. 
\begin{example} 
\begin{eqnarray}
\nonumber \ta_1(m)&=&-m,\ \ 
\ta_2(m) = -m(m-1)/10, \ \  
\ta_3(m) = (m+1)m(m-1)/210, \\
\nonumber
\ta_4(m)&=&-(19m+6)(m+1)m(m-1)/46200, \\
\nonumber
\ta_5(m)&=&-(155m^2-572m-48)(m+1)m(m-1)/21021000, \\
\nonumber \\
\nonumber\tb_1(m) &=& m, \ \ 
\tb_2(m) = -m(m+1)/10, \ \ 
\tb_3(m) = -(m+1)m(m-1)/210, \\
\nonumber
\tb_4(m) &=& -(19m-6)(m+1)m(m-1)/46200, \\
\nonumber
\tb_5(m) &=& (155m^2+572m-48)(m+1)m(m-1)/21021000, \\
\nonumber \\
\nonumber
\tc_1(m) &=& 0 \quad
\tc_2(m) = 3m(m+1)/70,  \quad
\tc_3(m) = -(m+1)m/350, \\
\nonumber
\tc_4(m) &=& -9(m+2)(m+1)m(m-1)/200200, \\
\nonumber
\tc_5(m) &=& 3(m+2)(m+1)m(m-1)/3503500, \\
\nonumber
\tc_6(m) &=& -(207m^2+207m+50)(m+2)(m+1)m(m-1)/4526522000,\\
\nonumber
\tc_7(m) &=& 9(107m^2+107m+4)(m+2)(m+1)m(m-1)/348542194000.
\end{eqnarray}
\end{example}

\begin{remark}\label{rem:neg}
  \begin{enumerate}
  \item[{\rm(i)}]  We can extend the recursion \eqref{eq:recursion} to
    negative  $n$. Then by the symmetry we have
    $T_{-n-1}(x)=T_n(x)$. From this, we can deduce $\tb_j(m)=\ta_j(-m)$
    and $\tc_j(m)=\tc_j(-m-1)$.  
  \item[{\rm(ii)}]  As a polynomial in $m$,
    $\ta_{j+1}(m)$ is divisible by $\ta_{j}(m)$ for $j\le3$, but this
    does not hold in general as the case $j=4$ shows. Likewise, $\tc_j(m)$
    divides $\tc_{j+1}(m)$ for $2\le j\le5$ but not for $j=6$.
  \item[{\rm(iii)}] The fact that $\tc_1(m)=0$ was given in
    \cite[Th.1]{Taneda}.
  \end{enumerate}
\end{remark}

\section{Yablonskii-Vorob'ev polynomial 
modulo a prime}

Fix a prime number $p>3$ once and for all.  We first establish a
special case of Theorem~\ref{th:mp+n}, namely for $m=1$ and
$n=0$. Once having this, the general case will be proved rather easily.

\begin{proposition}\label{th:p}
We have
\[
T_p(x) \equiv x^{d_p} \bmod p.
\]
\end{proposition}

\noindent{\it Proof.} The key ingredient is again the determinant formula
\eqref{eq:detexp};
$$
T_p(x)=\mu_p\tau_p(x).$$
Noting that $(2k-1)!!$ is prime to $p$ if
$k<(p+1)/2$ and  is divisible by $p$ exactly once if $(p+1)/2\le k\le
p$, we find the  exact power of $p$ which divides
$\mu_p=\prod_{k=1}^p(2k-1)!!$ is  $p^{(p+1)/2}$. So, if we put
$\mu_p'=p^{-(p+1)/2}\mu_p$, we have $\mu_p'\in\Z$ and
\begin{equation}
  \label{eq:detexp-p}
 T_p(x)=\mu_p'p^{(p+1)/2}\tau_p(x).
\end{equation}
We first show that the polynomial $p^{(p+1)/2}\tau_p(x)$ is realized as
a determinant of a matrix with entries which have $p$-integral coefficients.
To state this, we develop some notation.  Let $\Zp$ denote the local ring 
$\{ b/a \in \Q \mid a,b\in\Z, (a,p) = 1 \}$ which contains $\Z$ as a subring.
The maximal ideal of $\Zp$ generated by $p$ is denoted by $\p$.
Set $\p[x] = \{ \sum_{j\ge0} r_j x^j \in \Zp[x] \mid r_j \in \p \}$. By
``$\bmod p$'' of an element in $\Zp[x]$, we mean its image in the quotient
ring $\Zp[x]/\p[x]\simeq\F_p[x]$, where $\F_p$ is the field of $p$ elements.

Recall the polynomial $h_k(x)$ was defined by the generating function
\eqref{eq:generating}. Expanding $(d/d\lambda)e^{x \lambda+ \lambda^3/3}
= (x + \lambda^2)e^{x \lambda+ \lambda^3/3}$ 
we obtain the recursion
\[
(k+1) h_{k+1}(x) = x h_k(x) + h_{k-2}(x)
\quad\mbox{ for } k \ge 2,
\]
with $h_0=1, h_1=x$, and $h_2 = x^2/2$.  Multiplying both sides by
$k!$ and setting $\barh_k(x) = k! h_k(x)$, we have
\[
\barh_{k+1}(x) = x \barh_k(x) + k(k-1) \barh_{k-2}(x)
\quad\mbox{ for } k \ge 2,
\]
with $\barh_0=1, \barh_1=x$, $\barh_2 = x^2$.  This implies
inductively that $\barh_k(x)$ is a monic polynomial of degree $k$ with
integral coefficients.  In particular, we have 
\begin{equation}
  \label{eq:hkp}
  h_k(x)\in\Zp[x]\ \ \mbox{if}\  k<p 
\ \ \mbox{and}\ \ ph_k(x)\in\Zp[x]\ \ \mbox{if}\  p\le k<2p.
\end{equation}

Now define a matrix $(a_{ij})_{1 \le i,j \le p}$ by
\[
a_{ij} = \left\{
\begin{array}{ll}
h_{j-2i+p+1} & \mbox{if } i > (p+1)/2, \\
p h_{j-2i+p+1} & \mbox{if } i \le (p+1)/2.
\end{array}
\right.
\]
Then by \eqref{eq:hkp} and \eqref{eq:schur}, we have $a_{ij}\in\Zp[x]$
and 
\begin{equation}\label{eq:tau=det}
p^{(p+1)/2} \tau_p(x)
= \det(a_{ij})_{1 \le i,j \le p}.
\end{equation}
To compute this determinant modulo $p$, it is convenient to consider 
instead a modified matrix $(c_{ij})_{1 \le i,j \le p}$ which is obtained
from $(a_{ij})$ by a suitable permutation of rows: namely set
\[
c_{ij} = 
\left\{
\begin{array}{ll}
a_{kj} & \mbox{if } i=2k-1, \\
a_{k+(p+1)/2,j} & \mbox{if } i=2k.
\end{array}
\right.
\]
The inversion number of this permutation is
$\sum_{i=1}^{(p-1)/2} i = (p^2-1)/8$ and so 
\begin{equation}\label{eq:det=det}
\det(a_{ij}) = (-1)^{(p^2-1)/8} \det(c_{ij}).
\end{equation}
The following lemma supplies enough information for computing $\det(c_{ij})$
modulo $p$.
\begin{lemma}\label{lemma}
\begin{enumerate}
\item[{\rm(i)}]
If $i>j$, then $c_{ij} \in \p[x]$.

\item[{\rm(ii)}]
If $i$ is odd, then
$c_{ii}\equiv -x^p \bmod p.$

\item[{\rm(iii)}]
If $i$ is even, then
$c_{ii} = 1.$
\end{enumerate}
\end{lemma}
{\it Proof of Lemma.}~ If $i=2k-1$, then $k\le (p+1)/2$ and 
 $c_{ij}=a_{kj} = p h_{j-2k+p+1} = p h_{p-(i-j)}$. By \eqref{eq:hkp},
this belongs to $\p[x]$ if $i>j$, while for $i=j$ this is equal to $ph_p(x)=
\barh_p(x)/(p-1)!\equiv -\barh_p(x)\bmod p$ by Wilson's lemma. By \eqref{eq:hk},
the coefficient of $x^{p-3i}$ in $\barh_p(x)$ is
$p!/(3^ii! (p-3i)!)$, which is in $\p$ for $i \ge 1$ and hence 
$\barh_p(x)\equiv x^p\bmod p$. If $i=2k$, then
$c_{ij} = a_{k+(p+1)/2,j} = h_{j-2k} = h_{j-i}$. This is $0$ if $i>j$ and
$1$ if $i=j$. \qed

 From (i) of the lemma, the matrix $(c_{ij})$ modulo $p$ 
is upper-triangular, the diagonal entries of which are given by (ii) and (iii)
of the lemma. We therefore have 
\[
\det(c_{ij}) \equiv (-1)^{(p+1)/2}x^{p(p+1)/2} \bmod p.
\]
Combining this with \eqref{eq:det=det}, \eqref{eq:tau=det}, and
\eqref{eq:detexp-p}, we have 
$$ T_p(x)\equiv (-1)^{(p^2-1)/8+(p+1)/2}\mu'_p\,x^{d_p} \bmod p. $$
But we know that $T_p(x)$ is a monic polynomial of degree $d_p$, hence
the constant on the right should be congruent to $1$ modulo $p$ and we 
obtain the proposition. \qed
\begin{corollary}\label{th:cor-p}
We have $T_{p+1} \equiv x^{d_{p+1}} \bmod p$ and 
$T_{p-1} \equiv x^{d_{p-1}} \bmod p$.
\end{corollary}
{\it Proof.}~ From Proposition~\ref{th:p} we have $T'_p(x) \equiv 0
\bmod p$  since $d_p \equiv 0 \bmod p$. Thus  the recursion
\eqref{eq:recursion} reduces modulo $p$ to $T_{p+1} T_{p-1} \equiv x
T_p^2 \equiv x^{2d_p + 1}$.  Since $T_{p+1}(x)$ and $T_{p-1}(x)$ are
monic of degrees $d_{p+1}$ and $d_{p-1}$ respectively, and
$d_{p+1}+d_{p-1}=2d_p+1$, we get the formulas in the corollary.\qed

\noindent{\it Proof of Theorem~\ref{th:mp+n}.}~ Set
$S_n=x^{-(d_{n+p}-d_n)} T_{n+p}\bmod p$.  We know $S_0=1$ and $S_1=x$
by Proposition~\ref{th:p} and Corollary~\ref{th:cor-p}.  Noting that
$S'_n = x^{-(d_{n+p}-d_n)} T'_{n+p}$ and $S''_n = x^{-(d_{n+p}-d_n)} 
T''_{n+p}$ since $d_{n+p}-d_n =p(p+1+2n)/2
\equiv 0 \bmod p$, we have the same
recursion \eqref{eq:recursion} for $\{S_n\}_{n \ge 0}$.  Thus
we conclude $S_n\equiv T_n\bmod p$ for all $n$. Applying this inductively,
we establish Theorem~\ref{th:mp+n}. \qed
\begin{corollary} We have
$T_{p-1-i} \equiv x^{d_{p-1-i}-d_i} T_i \bmod p$.
\end{corollary}
{\it Proof.}~ We use the relation $T_{-n-1}(x)=T_n(x)$ as indicated
in Remark~\ref{rem:neg}. Theorem~\ref{th:mp+n} also holds for negative
indices and we obtain
$$T_{p-1-i}(x)=T_{-p+i}(x)\equiv x^{d_{-p+i}-d_i}T_i(x)
= x^{d_{p-i}-d_i}T_i(x).$$ \qed

Finally, we briefly mention what happens in the case when $p=2$ and $3$.
\begin{remark}
  Consider the general recursion \eqref{eq:gen-rec} in
  Remark~\ref{remark:3} with $a\in\Z$. For $p=3$, it is easy to see
  (using the fact that $T_n(x)$ is  ``almost'' a polynomial in $x^3$)
  that $$T_n(x)\equiv (x-a)^{d_n}\bmod 3\quad\mbox{if}\ n\equiv
  0,2\bmod 3$$
  and $$T_n(x)\equiv x(x-a)^{d_n-1}\bmod 3\quad\mbox{if}\
  n\equiv 1\bmod 3.$$
  In contrast, it trivially holds that
  $T_n(x)\equiv x^{d_n}\bmod 2$ if $a$ is even, while for odd $a$,
  numerical computation suggests that no periodic pattern for
  $T_n(x)\bmod 2$ exists and that irreducible factors of arbitrary
  high degree occur as $n$ gets bigger.
\end{remark}



\address

\end{document}